\documentclass{amsart}

\newtheorem{theorem}{Theorem}
\newtheorem{lemma}{Lemma}

\newtheorem{remark}{Proposition}

\begin{document}
\title[Spectrum of complex Jacobi matrices]{On limit sets for the discrete spectrum of complex Jacobi
matrices }

\author{Iryna Egorova}
\address{Institute for Low Temperature Physics and Engineering\\47\\ Lenin ave.\\ Kharkiv \\ Ukraine}
\email{egorova@ilt.kharkov.ua}

\author{Leonid Golinskii}
%\address{Institute for Low Temperature Physics and Engineering\\47\\ Lenin ave.\\ Kharkiv \\ Ukraine}
\email{golinskii@ilt.kharkov.ua}

\thanks{This work is supported in part by
INTAS Research Network NeCCA 03-51-6637.
%{\it To appear in ...}
}

\begin{abstract}
The discrete spectrum of complex Jacobi matrices that are compact
perturbations of the discrete laplacian is under consideration.
The rate of stabilization for the matrix entries sharp in the
sense of order which provides finiteness of the discrete spectrum
is found. The Jacobi matrix with the discrete spectrum having the
only limit point is constructed. The results can be viewed as the
discrete analogs of the well known theorems by B.S. Pavlov about
Schr\"{o}dinger operators on the half line with a complex
potential.

Bibliography : 26 references.
\end{abstract}

\maketitle

 \section{Introduction.}\setcounter{section}{1}
\setcounter{equation}{0}

$ \ \  \ $     In early 1960th B.S. Pavlov studied the discrete
spectrum of Schr\"odinger operators on the half line
\begin{equation}
\label{1.1} l_h y = -y^{\prime\prime} + q(x)y,\qquad y^{\prime}(0)
- h y(0) = 0,
\end{equation}
with a complex valued potential $q$ and a complex parameter $h$,
and came to the following results [1,2].

\noindent {\bf Theorem P1.} {\it The number of eigenvalues of the
operator $l_h$ is finite as long as the infinitely differentiable
potential $q(x)$ satisfies $$ |q(x)|\leq C\exp(-\varepsilon
x^\frac{1}{2}), \qquad x>0 $$ for some $\varepsilon >0$.}

\noindent {\bf Theorem P2.} {\it Whatever the numbers
$0<\beta<\frac{1}{2}$ and $\lambda>0$ are given, there is a real
infinitely differentiable potential $q(x)$ subject for some
$\delta>0$ to the bound
$$|q(x)|\leq C\exp(-\delta x^\beta),\qquad x>0 $$
and a complex number $h$ such that operator $(\ref{1.1})$ has
infinitely many eigenvalues with the only accumulation point
$\lambda$.}

\noindent The main objective of the present paper is to establish
the discrete analogs of Pavlov's theorems.

Let

\begin{equation}
    \label{1.3}
    J=\left(\begin{array}{ccccc}
                b_0&c_0&&&\\
                a_0&b_1&c_1&&\\
                &a_1&b_2&c_2&\\
                &&\ddots&\ddots&\ddots\\
            \end{array}\right)
\end{equation}
be a semi-infinite Jacobi matrix with complex entries and let
\begin{equation}
\lim_{n\to \infty}a_n=\lim_{n\to \infty}c_n=\frac{1}{2},\qquad
\lim_{n\to \infty}b_n = 0. \label{1.4}
\end{equation}
The operator $J$ in the Hilbert space $\ell^2( \mathbb{Z}_+)$,
$\mathbb{Z}_+ = \{0,1,2,...\}$, generated by such matrix is
thereby a complex perturbation of the discrete laplacian $$
J_0:\quad a_n=c_n=\frac{1}{2},\quad b_n=0.$$ According to the well
known Weyl's theorem the spectrum $\sigma(J)$ of $J$ is $\sigma(J)
= [-1,1]\bigcup \sigma_d(J)$, where a denumerable set
$\sigma_d(J)$ of eigenvalues of finite algebraic multiplicity,
called in what follows the {\it discrete spectrum} of the operator
$J$, may have accumulation points only on the interval $[-1,1]$.
The main object under consideration is a limit set for the
discrete spectrum $E_J:=\{\sigma_d(J)\}'$ for the operators
(\ref{1.3})-(\ref{1.4}).

The following question arises naturally: how to describe the class
$\{E_J\}_J$ of closed point sets on $[-1,1]$, which can emerge as
limit sets for discrete spectra of Jacobi matrices
(\ref{1.3})-(\ref{1.4}) or for certain subclasses of such
matrices. Our results can be viewed as the first step along the
way of solving this general problem. They concern primarily the
case of empty limit sets (finite discrete spectrum).

Let us say that $J$ (\ref{1.3}) belongs to the class $
\mathcal{P}(\beta)$ if
\begin{equation}
\label{1.5} \left|a_n -\frac{1}{2}\right| + \left|c_n -
\frac{1}{2}\right| + \left|b_n\right| \leq C_1\,\exp(-C_2\,
n^\beta),\qquad 0<\beta<1 ,\quad  C_1,\, C_2\,>0.
\end{equation}
To formulate our first result recall that the {\it convergence
exponent} or {\it Taylor-Besicovitch index} of a closed point set
$F\subset [-1,1]$ is the value $$\tau(F):=\inf\{\varepsilon >0:\
\sum_{j=1}^\infty |l_j|^\varepsilon <\infty\},$$ where $\{l_j\}$
are the adjacent intervals of $F$.
\begin{theorem} \label{theor1} Let
$J\in\mathcal{P}(\beta)$ where $0<\beta<\frac{1}{2}$. Then $E_J$
is a closed point set of the Lebesgue measure zero and its
convergence exponent satisfies
\begin{equation} \label{1.6} \dim E_J \leq
\tau(E_J)\leq\frac{1-2\beta}{1-\beta}\,,
\end{equation} where $\dim E_J$ is the Hausdorff dimension of
$E_J$. Moreover, if $J\in \mathcal{P}(\frac{1}{2})$ then
$E_J=\emptyset$ i.e. the discrete spectrum is finite.
\end{theorem}
The statement about finiteness of the discrete spectrum was proved
in a bit more general form in a recent paper [3, theorem 4.6].

It turns out that the exponent $1/2$ in theorem 1 is sharp in the
following sense.
\begin{theorem}
\label{theor2} For arbitrary $\varepsilon > 0$ and $-1<\nu<1$
there exists an operator $J\in \mathcal{P}(\frac{1}{2} -
\varepsilon)$ such that its discrete spectrum $\sigma_d(J)$ is
infinite and, moreover, the set $E_J$ consists of a single point:
$E_J=\{\nu\}$.
\end{theorem}

It is worth pointing out that for real Jacobi matrices $J$ with
$b_n=\overline{b_n}$, $a_n=c_n>0$ the abovementioned problem  is
simple. As the discrete spectrum in this case is real, the only
points $\pm 1$ (each one or both) can serve as its limit points.
As is well known [4], the finiteness of the first moment $$ \sum_n
n \left\{\left|a_n - \frac{1}{2}\right| + \left|b_n\right|\right\}
< \infty $$ already guarantees the finiteness of the discrete
spectrum of $J$.

The present paper appeared as an attempt to comprehend the
remarkable Pavlov's paper [1,2], which contains a beautiful
synthesis of ideas from the operator theory and the function
theory in the unit disk and which in our opinion was not fully
appreciated by experts. The discrete model we are dealing with is
comparatively simple that makes it possible to eliminate some
inaccuracy and simplify the arguments by using the well developed
techniques of perturbation determinants [5,6] in Section 2.
Section 3 follows the line of a seminal paper by L. Carleson [7]
on zero sets of analytic functions in the unit disk which provides
an analytic background for (\ref{1.6}) and enables one to complete
the proof of Theorem 1. In Section 4 we give some preliminary
information about the scattering problem for real matrices, which
is a key tool in the proof of Theorem \ref{theor2}in Section 5.
For the reader's convenience some results on perturbation
determinants for complex Jacobi matrices are collected in
Appendix.

Let us note in conclusion that the theory of complex Jacobi
matrices has undergone a splash of activity lately in connection
with the theory of formal orthogonal polynomials, Pad\'{e}
approximations and $J$-continued fractions [8-12].

\setcounter{section}{1} \setcounter{equation}{0}
\setcounter{lemma}{0} \setcounter{theorem}{0}
\renewcommand{\thetheorem}{\thesection.\arabic{theorem}}
\renewcommand{\thelemma}{\thesection.\arabic{lemma}}
\renewcommand{\theremark}{\thesection.\arabic{remark}}
\section {Perturbation determinants and Jost solution.}

In a recent paper [6] R. Killip and B. Simon applied the
perturbation determinants technique to the spectral theory of real
symmetric Jacobi matrices. Note that the most important results
from [6, sect.2] remain valid for the complex matrices as well.
For readers convenience we provide some information about
perturbation determinants in Appendix.

Throughout the rest of the paper we assume the finiteness of the
first moment
\begin{equation}
\label{2.1} \sum_{k=0}^{\infty} (k+1) \left\{ \,\left|a_k -
\frac{1}{2}\right| + \left|b_k\right| + \left|c_k -
\frac{1}{2}\right|\,\right\}<\infty,
\end{equation}
although some results (in a weakened form) hold under less
restrictive assumptions.

Under condition (\ref{2.1}) the operator $J$ is a trace class
perturbation of $J_0$, i.e., $\Delta J=J-J_0 \in \mathcal{S}_1$,
and so for $\lambda\notin [-1,1]$ $$(J-\lambda)(J_0-\lambda)^{-1}=
I + \Delta J(J_0 - \lambda)^{-1}, $$ and thereby the perturbation
determinant [5, Chap. IV.3]: $$\Delta(z,J):=\det(J-\lambda)(J_0 -
\lambda)^{-1}; \quad \lambda(z) = \frac{z+z^{-1}}{2},$$ is well
defined (see, e.g., [5, chap. IV.3]).

The main properties of the perturbation determinant are displayed
in the following statement [5, p.217; 6, Theorem 2.5].
\begin{remark}
 Under assumption $(\ref{2.1})$ the function  $\Delta$ belongs to the disk algebra
$\mathcal{A}$ of functions analytic in the unit disk $\mathbb{D} =
\{z:\,|z|<1\}$ and continuous up to the unit circle.
 $\Delta(z_0)=0$ for $|z_0|<1$ if and only if
 $\lambda(z_0)\in\sigma_d(J)$, and the order of zero equals the
 algebraic multiplicity of the eigenvalue  $\lambda(z_0)$.
 \end{remark}

A key role in the theory is played by the following approximation
relation for perturbation determinants. Denote by
$$
    J_m = \left(
    \begin{array}{ccccc}
        b_0&c_0&&&\\
        a_0&b_1&c_1&&\\
        &\ddots&\ddots&\ddots&\\
        &&\ddots&\ddots&c_{m-1}\\
        &&&a_{m-1}&b_m
    \end{array}
    \right),\qquad
 J_{0,m} = \left(
    \begin{array}{ccccc}
        0&\frac{1}{2}&&&\\
        \frac{1}{2}&0&\frac{1}{2}&&\\
        &\ddots&\ddots&\ddots&\\
        &&\ddots&\ddots&\frac{1}{2}\\
        &&&\frac{1}{2}&0
    \end{array}
    \right)
 $$ the matrices of order $m+1$ obtained by taking the first $m+1$
rows and columns of $J$ и $J_0$, respectively. Then ([6, formula
(2.59)])
\begin{equation}
\label{2.2} \Delta(z,J) = \lim_{m\to \infty} \det\left(\frac{J_m
-\lambda}{J_{0,m}-\lambda}\right)
\end{equation}
uniformly on compact subsets of $\mathbb{D}$ (for the proof of
(\ref{2.2}) see Appendix).

Along with the initial matrices $J,J_0$ consider the {\it
associated} matrices $J^{(n)},J^{(n)}_0$, obtained by dropping the
first $n+1$ rows and columns
$$
 J^{(n)}=\left(\begin{array}{ccccc}
                b_{n+1}&c_{n+1}&&&\\
                a_{n+1}&b_{n+2}&c_{n+2}&&\\
                &\ddots&\ddots&\ddots&\\
            \end{array}\right),
            \qquad
             J^{(n)}_0=J_0,
            $$
and $J^{(-1)}=J.$ The matrices $J_m^{(n)}$ are defined in a
natural way. By expanding the determinant $\det(J_m-\lambda)$ over
the first row and dividing by $\det(J_{0,m}-\lambda)$ we obtain
\begin{equation}\label{2.22}
\det\left(\frac{J_m-\lambda}{J_{0,m}-\lambda}\right)=(b_0-\lambda)
\det\left(\frac{J_{m-1}^{(0)}-\lambda}{J_{0,m-1}-\lambda}\right
)\det\left(\frac{J_{0,m-1}-\lambda}{J_{0,m}-\lambda}\right) -$$
$$- a_0c_0
\det\left(\frac{J_{m-2}^{(1)}-\lambda}{J_{0,m-2}-\lambda}\right
)\det\left(\frac{J_{0,m-2}-\lambda}{J_{0,m}-\lambda}\right).\end{equation}
Next, a direct calculation (cf. [6, formula (2.10)] shows that
$$ \lim_{m\to \infty}\det
\left(\frac{J_{0,m-j}-\lambda}{J_{0,m}-\lambda}\right)=(-2z)^j,\quad
j\in \mathbb{N}=\{1,2,\ldots\},$$ and so taking a limit in
(\ref{2.22}) as $m\to\infty$ with (\ref{2.2}) in mind leads to
\begin{equation}\label{2.3}\Delta(z,J)=(\lambda-b_0)(2z\Delta(z,J^{(0)}))-a_0c_0(4z^2\Delta(z,J^{(1)})).
\end{equation}
Since (\ref{2.3}) holds for any matrix $J$, write it for
$J=J^{(n)}$:
\begin{equation}\label{2.4}\Delta(z,J^{(n)})=(\lambda-b_{n+1})(2z\Delta(z,J^{(n+1)}))-
a_{n+1}c_{n+1}(4z^2\Delta(z,J^{(n+2)})).
\end{equation}
 Multiply (\ref{2.4}) by  $z^n$ and put
$\psi_n:=z^n\Delta(z,J^{(n)})$. After the shift of indices we have
\begin{equation}
\label{2.5} \psi_{m-1}(z) + 2b_m\psi_m(z) +4a_mc_m\psi_{m+1}(z) =
2\lambda\psi_m(z),\quad m\geq 0.
\end{equation}

Another important property of perturbation determinants is
contained in the following result [6, proposition 2.14], its proof
is given in Appendix.
\begin{remark}
 Under condition $(\ref{2.1})$ the relation
 \begin{equation}\lim_{m\to\infty}\Delta(z,J^{(m)})=1 \label{2.6}
 \end{equation}
holds uniformly on the closed unit disk $\overline{\mathbb{D}}$.
\end{remark}
Therefore for real symmetric $J$ the solution $\psi_m$ of
(\ref{2.5}) up to a constant factor (parameter) agrees with the
Jost solution of the equation $(J-\lambda)\psi = 0$ (see [13]).

In assumption (\ref{2.1}) denote
\begin{equation} \label{2.30}
 H(n):=\sum_{j=n}^\infty \left(|2b_m| + |4a_mc_m - 1|\right), \quad
H(n,m):=\prod_{j=n+1}^{n+m-1}(1+H(j)).
\end{equation}
It is easy to see that $\{H(n)\}\in \ell^1$ and the sequences
$H(\cdot)$ and $H(\cdot,m)$ decrease monotonically.

We can state the main result of the section as follows.

\noindent {\bf Theorem  2.3.}  {\it Under $(\ref{2.1})$ the Taylor
coefficients of the perturbation determinant
$$\Delta(z,J)=\sum_{j=0}^\infty \delta_j z^j$$ admit the bound
\begin{equation}
\label{2.7} |\delta_{j}| \leq \prod_{j=1}^{\infty}(1+H(j))
\,\sum_{m=\left[\frac{j}{2}\right]}^\infty \{|2b_m| + |4a_mc_m
-1|\},
\end{equation}
where $[x]$ is an integer part of $x$. In particular, $\Delta$
belongs to the space $W_+$ of absolutely convergent Taylor
series.}
%\end{theorem}

\noindent {\it Proof.} As a first step we show that $\psi_n$ obeys
the $\,$ "discrete integral"$\,$ equation of Volterra type. To
that end we introduce a function
 $$G(n,m,z)=\left\{ \begin{array}{cc}
2\,\frac{z^{m-n}-z^{n-m}}{z-z^{-1}},&m\geq n,\\0,&m<n,
\end{array}\right.$$
for which
\begin{equation}\frac{1}{2}\,G(n,m+1)+\frac{1}{2}\,G(n,m-1)-\lambda\,
G(n,m)=\delta(n,m) \label{2.8}
\end{equation}
holds with $\delta(n,m)$ the Kronecker symbol. Now let us multiply
(\ref{2.8}) by $\psi_m$, (\ref{2.5}) by $\frac{1}{2}G(n,m)$,
subtract the latter from the former and sum up the equalities
obtained over $m$ from $n$ to $N$. We have $$ \psi_n =
\sum_{m=n}^N\,\left\{-b_m G(n,m)+\left(\frac{1}{2} -
2a_{m-1}c_{m-1}\right)G(n,m-1)\right\}\psi_m +$$ $$+ \frac{1}{2}
G(n,N+1)\psi_N - 2a_Nc_NG(n,N)\psi_{N+1}. $$ By tending
$N\to\infty$ and taking into account the explicit form of the
kernel $G$ and (\ref{2.6}) we see that
\begin{equation}
\label{2.9} \psi_n(z)=z^n +\sum_{m=n+1}^\infty
M(n,m,z)\psi_m(z),\quad n=-1, 0,1,...,
\end{equation}
where $$M(n,m,z)=-b_m G(n,m,z)+\left(\frac{1}{2} -
2a_{m-1}c_{m-1}\right)G(n,m-1,z),$$ and $M(n,n,z)=0$. It is
convenient to transform (\ref{2.9}) by bringing in new variables
$$ \Delta(z,J^{(n)}) = z^{-n}\psi_n,\quad\hat
M(n,m,z)=M(n,m,z)z^{m-n},$$ and so $$
 \Delta(z,J^{(n)})=1 + \sum_{m=n+1}^\infty
\hat M(n,m,z)\Delta(z,J^{(m)}),$$ or otherwise,
\begin{equation}
\label{2.10} \Delta(z,J^{(n)})-1 = \sum_{m=n+1}^\infty \hat
M(n,m,z) + \sum_{m=n+1}^\infty \hat
M(n,m,z)\left(\Delta(z,J^{(m)})-1\right).
\end{equation}
It is clear that $\hat M$ is a polynomial on $z$ , and as for
$|z|\leq 1$
$$|G(n,m,z)z^{m-n}|=2|z^{2(m-n)} -
1||z-z^{-1}|^{-1}\leq2|z||m-n|,$$ $$
|G(n,m-1,z)z^{m-n}|\leq2|z|^2|m-n-1|<2|z||m-n|,$$ then
\begin{equation}
\label{2.11} |\hat
M(n,m,z)|\leq|z||m-n|\{|2b_m|+|1-4a_{m-1}c_{m-1}|\}, \quad |z|\leq
1.
\end{equation}
Condition (\ref{2.1}) ensures the series $$ g(n,z)=
\sum_{m=n+1}^\infty \hat M(n,m,z)$$ to converge absolutely and
uniformly inside $\overline{\mathbb{D}}$ and by (\ref{2.11}) $$
|g(n,z)|\leq\sum_{m=n+1}^\infty mh_m,\quad
h_m=|2b_m|+|1-4a_{m-1}c_{m-1}|.$$ The standard method of succesive
approximations (cf., e.g., [14, lemma 7.8]) gives $$
\left|\Delta(z,J^{(n)})-1\right|\leq\exp\left(\sum_{m=1}^\infty
mh_m\right)\cdot \sum_{m=n+1}^\infty mh_m,\quad n\geq -1,\ |z|\leq
1.$$

Consider now the Taylor series expansion of the function
$\Delta(z,J^{(n)})$:
\begin{equation}\label{2.12}\Delta(z,J^{(n)})=1+\sum_{j=1}^\infty\kappa(n,j)z^j
\end{equation}
and estimate its Taylor coefficients. By the Cauchy inequality
$$|\kappa(n,j)|\leq C\sum_{m=n+1}^\infty mh_m,$$ where by $C$ we denote
some positive constants which do not depend on the space variables
and the spectral parameter. In particular, $\kappa(n,j)\to 0$ as
$n\to\infty$ and fixed $j$.

The more accurate bound takes into account the dependence on the
second index $j$. If we plug (\ref{2.12}) in (\ref{2.4}) and match
the corresponding coefficients we have
\begin{equation}\label{2.13}
\kappa(n,j+1)=\kappa(n+1,j-1) - \sum_{m=n+1}^\infty \left
\{2b_m\kappa(m,j)+\left(4a_mc_m - 1\right)\kappa(m+1,j-1)\right\}
\end{equation}
for $j\geq2$,
$$ \kappa(n,1)=-2\sum_{m=n+1}^\infty b_m,\quad
\kappa(n,2)=-\sum_{m=n+1}^\infty \left\{2b_m \kappa(m,1) +(4a_mc_m
- 1)\right\}.$$ The induction on $j$ coupled with (\ref{2.13})
leads to
$$|\kappa(n,j)|\leq
H(n,j)H\left(n+1+\left[\frac{j}{2}\right]\right),\quad n\geq -1,$$
where $H(n)$, $H(n,m)$ are taken from (\ref{2.30}). Hence
\begin{equation}\label{2.14}
|\kappa(n,j)|\leq \prod_{j=1}^{\infty}(1+H(j))
\,\sum_{m=n+1+\left[\frac{j}{2}\right]}^\infty \{|2b_m| + |4a_mc_m
-1|\}.
\end{equation}
The desired bound (\ref{2.7}) is (\ref{2.14}) with $n=-1$. The proof is complete.% $\square$ .

\noindent {\bf Corollary 2.4.} {\it Let for $J$ $(\ref{1.3})$ the
moment of the order $n+1$ is finite
\begin{equation}
\label{2.15} M_{n+1}:=\sum_{k=0}^{\infty} (k+1)^{n+1} \left\{
\,\left|a_k - \frac{1}{2}\right| + \left|b_k\right| + \left|c_k -
\frac{1}{2}\right|\,\right\}<\infty.
\end{equation}
Then  the $n$'s derivative $\Delta^{(n)}(z)$ belongs to $W_+$ and
\begin{equation} \label{2.16}
\max_{z\in\overline{\mathbb{D}}} \left|\Delta^{(n)}(z)\right| \leq
C(J)\,\frac{4^n}{n+1}\,M_{n+1},
\end{equation}
where a positive constant $C(J)$ depends only on $J$.}

\noindent {\it Proof.} The statement is a simple consequence of
(\ref{2.15}), the series expansion
$$ \Delta^{(n)}(z)=\sum_{j=0}^\infty (j+1)\ldots
(j+n)\delta_{j+n}z^j $$ and bounds (\ref{2.7}). \hfill  $\square$

 \setcounter{section}{2} \setcounter{equation}{0}
\setcounter{lemma}{0} \setcounter{theorem}{0}

\section {Zero sets for classes of analytic functions in the unit
disk.}

Let $X$ be a class of functions from the disk algebra
$\mathcal{A}$. A closed point set $E$ on the unit circle
$\mathbb{T}$ is called the {\it set of uniqueness} for $X$, if
there is no nontrivial function $f\in X$, which vanishes on $E$.
Otherwise $E$ is a {\it zero set} for $X$. When $X=\mathcal{A}$,
then according to Fatou's theorem [15], a set $E$ is a zero set
for $ \mathcal{A}$ if and only if $E$ has the Lebesgue measure
zero: $|E|=0$. If $X$ is a subclass of $ \mathcal{A}$, the zero
sets may possess some additional properties. The investigation of
(metric) properties of zero sets goes back to A. Beurling [16] and
L. Carleson [7] (see also [17] for the extension of these
results). For instance let $\mathcal{A}_\infty$ be a set of all
functions from $\mathcal{A}$ which have all their derivatives in
the same class $ \mathcal{A}$. A closed point set
$E=\mathbb{T}\backslash \cup_j l_j$ with the adjacent arcs $l_j$
is the zero set for $ \mathcal{A}_\infty$ if and only if
\begin{equation}\label{3.1}
\sum_j |l_j|\log \frac{1}{|l_j|} < \infty \Leftrightarrow
\int_{\mathbb{T}} \left|\log\rho(\zeta, E)\right| dm<\infty, \quad
\rho(\zeta,E)=\mbox{dist}(\zeta, E),
\end{equation}
where $dm$ is the normalized Lebesgue measure on $\mathbb{T}$. We
are interested in certain subclasses of the class
$\mathcal{A}_\infty$, and for zero sets some more stringent
conditions than (\ref{3.1})  emerges.

 Let $g\in  \mathcal{A}_\infty$ и
\begin{equation}\label{3.2}
G_n(g)=G_n:=\max_{z\in\overline{\mathbb{D}}} |g^{(n)}(z)|,\quad
n\in \mathbb{Z}_+.
\end{equation}

The following result (Taylor's formula with the reminder term) can
be easily proved by induction on $n$.
\begin{lemma}
\label{lem3.1}
 Let $g\in  \mathcal{A}_\infty$. Then for an arbitrary $n\in \mathbb{N}$ and points $z,w$
in $\overline{\mathbb{D}}$ the inequality holds
\begin{equation}\label{3.3}
\left|g(z) -
\sum_{k=0}^{n-1}\frac{g^{(k)}(w)}{k!}(z-w)^k\right|\leq\frac{G_n}{n!}\left|z-w\right|^n.
\end{equation}
\end{lemma}

Consider the zero set of a function $g$ on $\mathbb{T}$:
$$F(g)=F:=\{\zeta\in\mathbb{T}\ :\ g^{(n)}(\zeta)=0\, \mbox{for
all}\, n\in\mathbb{Z}_+\}.$$
 Clearly $F$ is a closed point set of
measure zero which can be identified with a point set on
$[0,2\pi)$. Let $F_t=\{x:\mbox{dist}(x,F)\leq t\}$ and put
$\phi_F(t):=\phi(t)=|F_t|.$ The properties of $\phi$ as $t\to 0$
play a crucial role in the study of zero sets in some classes of
functions in $ \mathcal{A}_\infty$ with certain bounds on their
derivatives. The function $\phi$ is also known as the distribution
function of $d_F(x)=\mbox{dist} (x,F),\ x\in \mathbb{R}$. By the
change of variables formula for any measurable function $h$
\begin{equation}
\label{3.4} \int_{F_s} h(d_F(t))dt = \int_0^s h(u)d\phi(u)
\end{equation}
holds. For a function $g\in\mathcal{A}_\infty$ we define the value
(cf. [1])
\begin{equation}
\label{3.5} T(s):=\inf_{k\geq 0}\frac{G_k(g)}{k!}s^k,\quad s>0.
\end{equation}
\begin{lemma}
\label{lem3.2} Let $g\in\mathcal{A}_\infty$ and $g\neq 0$. Then $$
\int_0^{2\pi} \log T(s)d\phi(s)>-\infty.$$\end{lemma} {\it Proof.}
Put in (\ref{3.3}) $w=\exp(i\theta_0)$, $\theta_0 \in F$ and
$z=\exp(i\theta)$ . Then $$|g(e^{i\theta})|\leq \frac{G_n}{n!}
|\theta - \theta_0|^n.$$ The right hand side of this inequality
contains two parameters ($n$ и $\theta_0$). Taking the minimum we
have $$|g(e^{i\theta})|\leq \inf_{n\geq 0}\ \frac{G_n}{n!}\,
(d_F(\theta))^n=T(d_F(\theta)),\quad \log|g(e^{i\theta})|\leq \log
T(d_F(\theta)).$$
 It remains only to apply (\ref{3.4}) and the well-
known property $$\int_0^{2\pi}\log|g(e^{i\theta})|d\theta
>-\infty. $$
The lemma is proved.

Let us now introduce the main object related to the zeros of $g$
inside the disk: $$ E(g)=E:=\{\zeta \in \mathbb{T}:\
\exists\{z_n\}\in \mathbb{D},\ z_n\to\zeta,\ g(z_n)=0\}.$$
\begin{lemma}
\label{lem3.3} $E\subset F$.
\end{lemma}
{\it Proof.} Assume that $g^{(k)}(w)=0$, $k=0,1,...,n-2$ and
$g^{(n-1)}(w)\neq 0$ for some $w\in E$. By (\ref{3.3}) $$|g(z)
-\frac{ g^{(n-1)}(w)}{(n-1)!}(z-w)^{n-1}|\leq \frac{G_n}{n!}
|z-w|^n.$$ Putting $z=z_k$ and tending $k\to\infty$ leads to the
contradiction. \hfill  $\square$

\noindent {\it Proof of theorem \ref{theor1}} rests on Corollary
2.4. For matrices $J\in\mathcal{P}(\beta)$ (\ref{1.5}) condition
(\ref{2.15}) holds a fortiori for all $n$, and it is easy to find
 bounds for the moments $M_r$. As
$$ \left|a_k - \frac{1}{2}\right| + \left|b_k\right| + \left|c_k -
\frac{1}{2}\right|\leq C\exp(-C(k+1)^\beta),\quad 0<\beta<1,$$ (by
$C$ we denote positive constants which depend only on the original
matrix $J$), we have
$$M_r\leq C\sum_{k=0}^\infty(k+1)^r
\exp\left(-\frac{C}{2}(k+1)^\beta\right)\exp\left(-\frac{C}{2}(k+1)^\beta\right).$$
An undergraduate analysis of the function
$u(x)=x^r\exp(-\frac{C}{2}x^\beta)$ gives $$\max_{x\geq 0
}u(x)=u(x_0)=\left(\frac{2r}{C\beta}\right)^\frac{r}{\beta}e^{-\frac{r}{\beta}}=
\left(\frac{2}{C\beta e}\right)^\frac{r}{\beta} r
^\frac{r}{\beta},\quad
x_0=\left(\frac{2r}{C\beta}\right)^\frac{1}{\beta},$$ so that
$$ M_r\leq B\left(\frac{2}{C\beta e}\right)^{\frac{r}{\beta}}
r^{\frac{r}{\beta}},\qquad
B=C\sum_{k=0}^\infty\exp\left\{-\frac{C}{2}(k+1)^\beta\right\}.$$
Hence $$G_n(\Delta)=\max_{z\in\overline{\mathbb{D}}}
\left|\Delta^{(n)}(z)\right|\leq
C\frac{4^n}{n+1}\left(\frac{2}{C\beta
e}\right)^{\frac{n+1}{\beta}}(n+1)^{\frac{n+1}{\beta}},\quad n\geq
0.$$ The inequalities
$$\left(\frac{n+1}{n}\right)^\frac{n}{\beta}<e^\frac{1}{\beta},\qquad
(n+1)^\frac{1}{\beta}<e^\frac{n+1}{\beta},\qquad
(n+1)^\frac{n+1}{\beta}<e^\frac{1}{\beta}n^\frac{n}{\beta}e^\frac{n+1}{\beta}$$
lead to  $$G_n(\Delta)\leq C C_1^n n^\frac{n}{\beta},\quad n\geq
0.$$ In other words, the perturbation determinant for
$J\in\mathcal{P}(\beta)$ belongs to the Gevre class
$\mathcal{G}_\beta$.

Let us go over to the function $T$ (\ref{3.5}). As $n!>n^ne^{-n}$,
then $$T(s)\leq C \inf_n (C_1 e s )^n n^{\frac{n}{\beta} -
n}=C\inf_n t^n n^{\alpha n},$$ where $t=C_1 e s <1/2$ for small
enough $s$, $\alpha = \beta^{-1} - 1 > 0$. An elementary analysis
of the function $v(x)=t^x x^{\alpha x}$ shows that
\begin{equation}\label{3.20} \min_{x\geq 0}
v(x)=v(x_1)=\exp\left(-\frac{\alpha}{e}\
t^{-\frac{1}{\alpha}}\right),\qquad x_1=\frac{1}{e}\
t^{-\frac{1}{\alpha}} \gg 1.\end{equation} Although the number
$x_1$ is not in general an integer, by putting $n=[x_1]$, it is
not hard to make sure that $$T(s)\leq
C\exp\left\{-C\left(\frac{1}{s}\right)^\frac{\beta}{1-\beta}\right\},\qquad
s\leq s_0.$$ Applying lemma \ref{lem3.2}, we conclude that
\begin{equation}\label{3.6}
\int_0 \frac{d\phi_F(s)}{s^\frac{\beta}{1-\beta}} < \infty, \qquad
F=F(\Delta).
\end{equation}
Suppose first that $0<\beta<1/2$. It is proved in [7] that the
convergence of integral (\ref{3.6}) is equivalent to convergence
of the series $$\sum_{j=1}^\infty |l_j|^\frac{1-2\beta}{1-\beta}
<\infty,$$ where $\{l_j\}$ are adjacent intervals of the closed
point set $F$ (as a set on $[0,2\pi)$), and thereby
$$\tau(F)\leq\frac{1-2\beta}{1-\beta}.$$ The inequality $\dim
F\leq\tau(F)$ follows from the general theory of fractal dimension
(see [18]). By lemma \ref{lem3.3} $E=E(\Delta)\subset F$, and so
due to monotonicity of dimensions $$\dim
E\leq\tau(E)\leq\frac{1-2\beta}{1-\beta}\,.$$ Inequality
(\ref{1.6}) stems out of an obvious relation between $E_J$ и $E$.

Consider now the case $\beta=\frac{1}{2}$, so that
\begin{equation}\label{3.7}
\int_0\frac{d\phi_F(s)}{s} <\infty.\end{equation} We show that the
latter is impossible for nonempty $F$. Suppose that the adjacent
arcs $l_j=(\alpha_j,\beta_j)$ are labelled in the order of
decreasing length, i.e., $|l_1|\geq|l_2|\geq...$. Let
$0<t_1<t_2<\frac{1}{2}|l_1|=\frac{\beta_1 - \alpha_1}{2}$ . By the
definition of $\phi$ $$\phi(t_2) - \phi(t_1) = |\{x:\
t_1<\mbox{dist}(x,F)\leq t_2\}|.$$ For the interval
$I(t_1,t_2)=(\alpha_1 + t_1, \alpha_1 +t_2)$ the relation
$$I(t_1,t_2)\subset\{x:\ t_1<\mbox{dist}(x,F)\leq t_2\}$$ takes
place, and so $\phi(t_2) - \phi(t_1)\geq t_2 - t_1$. Hence the
measure $d\phi_F$ dominates the Lebesgue measure and integral
(\ref{3.7}) diverges. So we arrive at the following conclusion: if
$J\in\mathcal{P}(\frac{1}{2})$, then $F$ (as well as $E$) is
empty, which exactly means finiteness of the discrete spectrum
$\sigma_d(J)$. The proof is complete.

\noindent {\bf Remark 1}. Corollary 2.4 enables one to obtain some
conditions on the  "size" of the set $E_J$ in much more general
situations. For instance, if $J$ satisfies (\ref{2.15}), then
$\Delta\in\mathcal{A}_\infty$, and so for the set $E(\Delta)$
(\ref{3.1}) is true, and the same condition (taken over to the
interval $[-2,2]$) holds for $E_J$. As another example, let us
take the class of Jacobi matrices $J$ subject to $$ \left|a_k -
\frac{1}{2}\right| + \left|b_k\right| + \left|c_k -
\frac{1}{2}\right|\leq
C\exp\left\{-C(\log(k+1))^\gamma\right\},\qquad \gamma>1,$$ which
contains the class $\mathcal{P}(\beta)$. It is not hard to verify
that for perturbation determinants of such matrices $$
\max_{z\in\overline{\mathbb{D}}} \left|\Delta^{(n)}(z)\right|\leq
C\exp\{Cn^p\}, \qquad p=\frac {\gamma}{\gamma-1}>1, $$ holds. The
functions with such kind of bound for their derivatives were
studied by B. Taylor and D. Williams in [19] , wherein the
condition on the zero sets $$ \int_{\mathbb{T}}
\left|\log\rho(\zeta, E)\right|^q dm<\infty, \quad
\frac1{q}+\frac1{p}=1, $$ was established.

On the other hand, the finiteness of the discrete spectrum for
$J\in\mathcal{P}(\frac12)$ is caused by quasi-analyticity of the
Gevre class $\mathcal{G}_{1/2}$, i.e., the lack of nontrivial
function $f$, such that $f^{(n)}(\zeta_0)=0$ for all $n\geq0$ and
some $\zeta_0\in\mathbb{T}$. The well-known criterion of
quasi-analyticity [20] makes it possible to prove slightly more
general than (\ref{1.5}) with $\beta=1/2$ results concerning the
finiteness of the discrete spectrum.

\noindent {\bf Remark 2}. The perturbation determinant of
$J\in\mathcal{P}(\frac12)$ may have boundary zeros  on the unit
circle, that can be observed on a simple example $$ J:\qquad
a_n=c_n=\frac{1}{2},\quad b_0=e^{i\theta}, \ \ b_1=b_2=\ldots=0.$$
It is easy to compute $\Delta(z,J)=1-e^{i\theta}z$. These points
(more precisely, their images under the Zhukovsky transformation)
are known as the spectral singularities. In the above example with
$\theta\not=0, \pi$ the matrix $J$ has empty discrete spectrum and
one spectral singularity at $\cos\theta\in (-1,1)$. Under
condition (\ref{2.1}) the spectral singularities form a closed set
of measure zero. It is clear that $J\in\mathcal{P}(\frac12)$ may
have at most finite number of spectral singularities (cf. [3]).
Indeed, otherwise we could find a point $\zeta_0\in\mathbb{T}$,
for which $\Delta^{(n)}(\zeta_0,J)=0$ for all $n\geq0$ (cf.
(\ref{3.3})), that contradicts to quasi-analyticity of
$\mathcal{G}_{1/2}$.

\setcounter{section}{3} \setcounter{equation}{0}
\setcounter{lemma}{0} \setcounter{theorem}{0}

\section {On scattering data for real symmetric Jacobi matrix.}
In this section we study the relation between Fourier coefficients
of the scattering function and the rate of stabilization of the
Jacobi matrix entries in the case when both the discrete spectrum
and resonance are absent. Our consideration is based  on the
scattering problem for semi-infinite real symmetric Jacobi matrix
[4, 13, 21]. The main result of the section is lemma \ref{lem4.2}.

Let
\begin{equation}
    \label{4.1}
    J=\left(\begin{array}{ccccc}
                b_1&a_1&&&\\
                a_1&b_2&a_2&&\\
                &\ddots&\ddots&\ddots&\\
                &&\ddots&\ddots&\ddots\\
            \end{array}\right)
\end{equation}
be a Jacobi matrix with the entries $a_n>0$, $b_n=\overline{b_n},$
subject to (\ref{1.6}). Put $a_0=1.$ It is well known (see [4, 13,
21]), that equation
\begin{equation}
\label{4.2} a_{n-1}y_{n-1} + b_ny_n +a_ny_{n+1}
=\frac{1}{2}(z+z^{-1})y_n,\quad z\in\overline{\mathbb{D}},\ n\geq
1
\end{equation}
has a solution $f_n(z)$, called the {\it Jost solution}, with the
following properties:

\noindent (i)$\quad \lim_{n\to\infty} f_n(z)z^{-n} =1$ uniformly
on $z\in\overline{\mathbb{D}}$ ;

\noindent (ii) for all $n\geq 0$ the function $f_n(z)$ belongs to
the algebra $W_+$ of absolutely convergent Taylor series, and
$f_n( \overline{\zeta})=\overline{f_n(\zeta)}$,
$\zeta\in\mathbb{T}$;

\noindent (iii) a function $f_0(z)$, known as the {\it Jost
function} has at most finite number of zeros in
$\overline{\mathbb{D}}\setminus \{\pm 1\}$, all of them are real
and simple, and their images in the $\lambda$-plane are points of
the discrete spectrum of $J$;

\noindent(iv) two solutions $f_n(\zeta)$ and $
\overline{f_n(\zeta)}$ of (\ref{4.2}) are linearly independent for
$\zeta \in \mathbb{T}\setminus\{\pm 1\}$ and
\begin{equation}\label{4.3}
\langle f, \overline{f}\rangle = \frac{\zeta^{-1} - \zeta}{2},
\end{equation}
where $\langle f, \overline{f}\rangle : =
a_{n-1}(f_{n-1}\overline{f_n} -f_{n}\overline{f_{n-1}})$is the
Wronskian of these solutions.

In view of (i) и (ii) the following representation holds
\begin{equation}
\label{4.4} f_n(z)=\sum_{m=n}^\infty K(n,m) z^m,\quad
z\in\overline{\mathbb{D}},\quad n\in\mathbb{Z}_+,
\end{equation}
where the entries $K(n,m)$ of the matrix $K$, called {\it the
transformation operator}, are real and obey the relations [14, 22]
\begin{equation} \label{4.5}
 a_n=\frac{K(n+1,n+1)}{2 K(n,n)},\ \mbox{that is}\
 K(n,n)=\prod_{j=n}^\infty(2 a_j)^{-1},
 \end{equation}
 \begin{equation}\label{4.6}
 b_n=\frac{K(n,n+1)}{2K(n,n)}-\frac{K(n-1,n)}{2K(n-1,n-1)},\quad
 n\geq1.
 \end{equation}
Note that the Jost function agrees (up to a constant factor which
does not depend on the spectral parameter) with the perturbation
determinant of $J$ (see, e.g., [6, formula (2.64)])
 $$f_0(z)=(\prod_{j=1}^\infty\, 2a_j)^{-1}\Delta(z,J).
 $$
Let $S(\zeta)= \overline{f_0(\zeta)}(f_0(\zeta))^{-1}$ be the so
called {\it scattering function} of the matrix $J$, $\zeta \in
\mathbb{T}$. As $f_0(\zeta)\neq 0$,
 $\zeta\in\mathbb{T}$ (under the assumption of the lack of resonance  $f_0(\pm1)\neq
0$), then $S(\zeta)$ is a continuous function on $\mathbb{T}$, and
moreover $ \overline{S(\zeta)}=S(\overline{\zeta})=S^{-1}(\zeta)$.
Put
\begin{equation}
\label{4.8} F(n):= -\frac{1}{2\pi}\int_{-\pi}^\pi
S(e^{i\theta})e^{in\theta} d\theta.
\end{equation}

In the case when there is no discrete spectrum, the Fourier
coefficients $F(n)$ and the matrix entries $K(n,m)$ of the
transformation operator are connected by the equation of the
inverse scattering problem (the Marchenko equation, see [4, 13])
\begin{equation}
\label{4.9} \frac{\delta(n,m)}{K(n,n)}=K(n,m) + \sum_{l=n}^\infty
K(n,l)F(l+m),\quad n,m\in\mathbb{Z}_+.
\end{equation}
As is known [14, formula (10.87)], under condition (\ref{1.6}) the
inequality
\begin{equation}
\label{4.10} \sum_{n=1}^\infty n|F(n+2) - F(n)|<\infty
\end{equation}
holds.

Consider the discrete analog of the Gelfand-Levitan operator $
\mathcal{F}_n :\ell^1(\mathbb{Z}_+) \rightarrow
\ell^1(\mathbb{Z}_+)$ $$ (\mathcal{F}_n y)_j = \sum_{m=0}^\infty
F(2n+m+j) y_m,\quad j\in\mathbb{Z}_+. $$ Here $n\geq 1$ is fixed.
\begin{lemma}
\label{lem4.1} Under assumption $(\ref{4.10})$ the operator $
\mathcal{I}+\mathcal{F}_n$ is invertible in the space
$\ell^1(\mathbb{Z}_+)$ for all $n\geq 1$ and

 $$\|(\mathcal{I}+\mathcal{F}_n)^{-1}\|\to 1 \ \  \mbox{as}\ \  n\to
\infty.$$
\end{lemma}
{\it Proof.} Put
\begin{equation}\label{4.14}
\hat F(n):=\sum_{j=0}^\infty|F(n+j) - F(n+j+2)| =
\sum_{k=n}^\infty |F(k) - F(k+2)|.\end{equation} Since $F(n)\to 0$
as $n\to \infty$, then $$ |F(k)|\leq \sum_{j=0}^\infty |F(k+2j) -
F(k+2j+2)|,$$ so that
\begin{equation}\label{4.18}
|F(k)| + |F(k+1)|\leq\hat F(k).
\end{equation} Hence $\{F(k)\}\in
\ell^1(\mathbb{Z}_+)$, i.e., the operator  $ \mathcal{F}_n$ is
compact.  According to the Fredholm alternative,
$\mathcal{I}+\mathcal{F}_n$ is invertible as soon as the equation
\begin{equation}\label{4.25}
\left(\mathcal{I}+\mathcal{F}_n\right)\,g = 0
\end{equation}
has only trivial solution. As
$S\left(e^{-i\theta}\right)=\overline{S\left(e^{i\theta}\right)}$,
we have $F(n)\in\mathbb{R}$, and it suffices to restrict ourselves
with the real solutions of (\ref{4.25}). Let $g=\{g(k)\}\in\ell^1$
be such solution. Take the function $$\widetilde
g(z):=\sum_{k=0}^\infty g(k)z^k.$$ Then
$$(g,g)=\sum_0^\infty |g(k)|^2 = \frac{1}{2\pi}\int_{-\pi}^\pi
\left|\widetilde g\left(e^{i\theta}\right)\right|^2\,d\theta,$$
$$(\mathcal{F}_n g,g) =-\frac{1}{2\pi}\int_{-\pi}^\pi\,
S\left(e^{i\theta}\right)e^{2in\theta}\,\widetilde
g^2\left(e^{i\theta}\right)\,d\theta$$ and so
\begin{equation}\label{4.26} 0=\left((\mathcal{I}+\mathcal{F}_n)g,g\right)=\frac{1}{2\pi}\int_{-\pi}^\pi\,
\left(1 - \Phi\left(e^{i\theta}\right)\right)\left|\widetilde
g\left(e^{i\theta}\right)\right|^2\,d\theta, \end{equation} where
$$\Phi\left(e^{i\theta}\right)=S\left(e^{i\theta}\right)e^{2in\theta}\,\widetilde
g^2\left(e^{i\theta}\right)\left|\widetilde
g\left(e^{i\theta}\right)\right|^{-2}.$$ Suppose that $\widetilde
g\neq 0$. Since $\left|\Phi\left(e^{i\theta}\right)\right|=1$ ,
then (\ref{4.26}) implies $\Phi\left(e^{i\theta}\right)\equiv 1$,
and so by the definition of the scattering function
$$\frac{\widetilde g^2(e^{i\theta})}{f_0^2(e^{i\theta})}\
e^{2in\theta} = \frac{\left|\widetilde
g(e^{i\theta})\right|^2}{\left|f_0(e^{i\theta})\right|^2}.$$
Thereby for the function $h_n=\widetilde
g^2\,f_0^{-2}\,e^{2in\theta}$ we have
$h_n(e^{i\theta})=|h_n(e^{i\theta})|$. Thanks to the lack of the
discrete spectrum and resonance, $f_0(z)\neq 0$ for
$z\in\overline{\mathbb{D}}$, and hence $h_n \in W_+$. Then by the
symmetry principle and the uniqueness theorem for analytic
functions  $h_n\equiv \mbox{const}$ in $\overline{\mathbb{D}}$.
But for $n\geq 1$ $h_n(0)=0$, that is, $\tilde g\equiv 0$. The
latter means that (\ref{4.25}) has only trivial solution and so
$\mathcal{I} + \mathcal{F}_n$ is invertible for each $n\geq 1$.
The second statement is a simple consequence of inequalities
$$\|\mathcal{F}_n\|\leq \sum_{j=0}^\infty |F(2n+j)| \rightarrow
0,\quad n\to\infty$$ and $$\|(\mathcal{I} + \mathcal{F}_n)^{-1} -
\mathcal I\|\leq\frac{\|\mathcal F_n\|}{1 - \|\mathcal F_n\|}.$$
The proof is complete. \hfill $\square$

Denote \begin{equation}\label{4.11}\tau(n,j) = K(n,n+j)K(n,n) -
\delta(n,n+j).\end{equation} Then (\ref{4.9}) takes the form
\begin{equation}
\label{4.12} \tau(n,j) + \sum_{m=0}^\infty F(2n + m +j)\tau(n,m) =
-F(2n +j),\quad j\geq 0,
\end{equation}
and hence \begin{equation}\label{4.13}\{\tau(n,j)\}_{j\geq 0} =
-(\mathcal{I}+\mathcal{F}_n)^{-1}\{F(2n+j)\}_{j\geq 0}
.\end{equation}

We need a slight modification of ([14, formula (10.110)]).
\begin{lemma}
\label{lem4.2}  Under assumption $(\ref{4.10})$ the bound holds
\begin{equation}\label{4.15}
|2a_n - 1| + |b_n|\leq C \{|F(2n-1) - F(2n+1)| + |F(2n) - F(2n+2)|
+ \hat F^2(2n-2)\},
\end{equation}
where $\hat F(n)$ is defined in $(\ref{4.14})$.
\end{lemma}
{\it Proof.} In the space $\ell^1(\mathbb Z_+)$ consider the
equation $(\mathcal{I} +\mathcal{F}_n) x = y$ or in other words
\begin{equation}\label{4.16}
x(j) +\sum_{m=0}^\infty F(2n+j+m)x(m) =y(j), \quad x=\{x(j)\},\
y=\{y(j)\},\ j\geq 0.
\end{equation}
Since $ \mathcal{I} +\mathcal{F}_n$ is invertible for $n\geq 1$
then the solution (\ref{4.16}) satisfies
\begin{equation} \label{4.17} \|x\|_1\leq K \|y\|_1,\quad
K=\sup_{n\geq 1}\|(\mathcal{I} +\mathcal{F}_n)^{-1}\|<\infty.
\end{equation}
By plugging (\ref{4.17}) in (\ref{4.16}) , we obtain the bound for
coordinates $$ |x(j)|\leq \sup_{k\geq 2n+j}|F(k)| \|x\|_1 +
|y(j)|\leq K \sup_{k\geq 2n+j}|F(k)|\|y\|_1 +|y(j)|.$$ By
(\ref{4.18}) and on account of  $\{\hat F(n)\}$ being
non-increasing, we come to inequality
\begin{equation}\label{4.19}
|x(j)|\leq K\hat F(2n+j)\|y\|_1 + |y(j)|.
\end{equation}

We  apply (\ref{4.19}) twice. First, to   equation (\ref{4.12})
with $$x(j)=\tau(n,j),\quad y(j)=-F(2n+j).$$ By
(\ref{4.10})-(\ref{4.11}) $$\|y_n\|_1\leq \sum_{k=0}^\infty \hat
F(k)=C<\infty,\quad |y_n(j)|\leq \hat F(2n+j), $$ and therefore,
(\ref{4.19}) takes the form
\begin{equation}\label{4.20}
|\tau(n,j)|\leq C \hat F(2n+j),\quad j\geq 0.\end{equation} Next,
write (\ref{4.12}) with indices $n$ and $n+1$ and subtract the
latter from the former: $$(\mathcal{I} +\mathcal{F}_n)(\tau_n -
\tau_{n+1}) = y_n,\quad \tau_n = \{\tau(n,j)\}_{j\geq 0},\quad
y_n=\{y_n(j)\}_{j\geq 0},$$ where
$$y_n(j)=F(2n+2+j) - F(2n+j) +\sum_{m=0}^\infty \{F(2n+j+m) -
F(2n+2+j+m)\}\tau(n+1,m).$$ It follows directly from (\ref{4.20})
that $$ |y_n(j)|\leq |F(2n+j) - F(2n+2+j)| + C\hat
F(2n)\hat(2n+j),\quad \|y_n\|\leq C\hat F(2n),
$$ and so (\ref{4.19}) implies
\begin{equation}\label{4.21} |\tau(n,j) - \tau(n+1,j)|\leq
|F(2n+j) - F(2n+2+j| + C \hat F(2n)\hat F(2n+1),\quad j\geq 0.
\end{equation}

The desired bound (\ref{4.15}) flows out of equalities
(\ref{4.5}), (\ref{4.6}) and definition (\ref{4.11}). Indeed, by
(\ref{4.5}) $$ 0< C_1\leq K^2(n,n)\leq C_2<\infty,\quad K(n,n)\to
1,\quad n\to\infty,$$ and so by (\ref{4.21}) with $j=0$ $$|2a_n -
1|=\frac{|K(n+1,n+1) -K(n,n)|}{K(n,n)}\leq C |\tau(n,0) -
\tau(n+1,0)|\leq $$\begin{equation}\leq C\{|F(2n) - F(2n+2)| +
\hat F^2(2n)\}. \label{4.22}
\end{equation}
Next, since $|\tau(n,1)|\leq C$ (see(\ref{4.20})), we have by
(\ref{4.6}) $$|2b_n| = \frac{|4a_{n-1}^2 \tau(n-1,1) -
\tau(n,1)|}{K^2(n,n)}\leq C\{|\tau(n-1,1) - \tau(n,1)| +
|4a_{n-1}^2 - 1|\}.$$ Taking into account (\ref{4.21}) for $j=1$,
(\ref{4.22}) and monotonicity of $\hat F(n)$ we see that
\begin{equation}\label{4.23}
|b_n|\leq C\{|F(2n-1) - F(2n+1)| +|F(2n) - F(2n+2)| + \hat
F(2n-2)\hat F(2n-1)\}.\end{equation} Finally, (\ref{4.22}) coupled
with (\ref{4.23}) gives (\ref{4.15}), that completes the proof of
lemma \ref{lem4.2}.

 \setcounter{section}{4} \setcounter{equation}{0}
\setcounter{lemma}{0} \setcounter{theorem}{0}
\section{Pavlov's example.}

The proof of Theorem \ref{theor2} is based on a specific
"oscillatory integral", invented by Pavlov in  [2, \S 3]
 \begin{equation}
 \label{5.1}
 V(z) = \int_0^z
 \exp\{-\chi(\xi)\}\cos(\gamma\chi(\xi))d\xi,\end{equation}
 where $\chi(\xi)=\frac{1}{32}(1 + \xi^2)^{\gamma - 1}$ and $\gamma$
 is a small enough parameter,
 $0<\gamma<1$. The following properties of such integrals are of particular
 interest:

 \noindent (i)$\quad V\in\mathcal{A}_\infty $ and is analytic at every point
 $\zeta\in\mathbb{T}\setminus \{\pm i\};$

 \noindent (ii)$\quad V(0) = 0$, $|V(z)|\leq 1$ for
 $z\in\overline{\mathbb{D}}$, and $V(z)\neq 0$ for $z\neq 0;$

 \noindent (iii)$\quad \mbox{Im}\,V(z)\ \mbox{Im}\,z>0$ for $z\notin
 \mathbb{R};$

 \noindent (iv) There is a sequence of points on the imaginary axis
 $z_k=i\Im z_k$, $0<\Im z_k<1$, such that $\Im z_k\uparrow 1$ and $V(z_k) = V(i)$. Moreover
 $V(\bar{z}_k) = V(-z_k) = \overline{V(i)}.$

We will turn to some other properties of Pavlov's function in due
course. Put
 $$ f(\lambda) = -\frac{1}{V(-z(\lambda))},\quad
 \lambda(z)=\frac{1}{2}\left(z+\frac{1}{z}\right).$$
 It is clear that $f$ is analytic in $\{\mbox{Im}\,\lambda>0\}$ and continuous
 (in fact infinitely differentiable) in $\{\mbox{Im}\,\lambda\geq
 0\}$ with
 $\mbox{Im}\,f(x) = 0$ for  $
 x\in(-\infty,-1]\cup[1,\infty)$ and $\mbox{Im}\,f(x)>0$ for
 $\lambda\in \{\mbox{Im}\,\lambda >0\}\cup(-1,1)$. Hence,
 \begin{equation}
 \label{5.2}
 f(\lambda)=\alpha\lambda + \beta +\int_{-1}^1
 \frac{\mbox{Im}\,f(x)}{x-\lambda}\,dx,\quad \alpha\geq 0,\
 \beta=\overline{\beta}.
 \end{equation}
 A coefficient $\alpha$ can be found from the limit relation
 $$\alpha =\lim_{y\to +\infty}\frac{f(i y)}{i
 y}= -\lim_{t\downarrow 0}\frac{2}{i(t^{-1} - t )V(it)},
 $$ and since $V(it)=V'(0)it +O(t^2)$ for $t\to 0$, then
 $$\alpha=\frac{2}{V'(0)}=\frac{2e^{1/32}}{\cos \gamma}>0.$$
 Define
 \begin{equation}\label{5.21} \widetilde{m}(\lambda):=A\int_{-1}^1
 \frac{\mbox{Im}\,f(x)}{x-\lambda}dx = A(f(\lambda)-\alpha\lambda
 -\beta),\quad A^{-1}:=\int_{-1}^1
 \mbox{Im}\,f(x)dx.\end{equation}
 The function $\widetilde{m}(\lambda)$ is the Weyl function for a  certain
 Jacobi matrix
$$ \tilde{J}=\left(\begin{array}{cccc}
                b_1&a_1&&\\
                a_1&b_2&a_2&\\
                &\ddots&\ddots&\ddots\\
            \end{array}\right),\quad a_j>0,\ b_j=\overline{b_j}.$$
Finally, put \begin{equation}\label{5.3} a_0^2:=-\frac{1}{\alpha
A},\quad b_0:=-\frac{\beta}{\alpha}-\frac{1}{\alpha
\overline{V(i)}}\,,\quad a_0,b_0\in \mathbb{C}\,,
\end{equation}
where $\alpha$ and $\beta$ are taken from representation
(\ref{5.2}), and construct a complex Jacobi matrix
$$J=\left(\begin{array}{cc} b_0&a_0\\ a_0&\tilde{J}\\
\end{array}\right).$$
Our goal is to prove that $J$ satisfies theorem \ref{theor2} with
$\nu = 0$. \footnote{At the final stage of the proof we show that
the general case only insignificantly differs from this particular
one.}

Let us begin with the discrete spectrum of $J$. Denote by
$R(\lambda)=(J-\lambda)^{-1}$ the resolvent of $J$ and let
$m(\lambda)=(R(\lambda)e_0,e_0)$ be its Weyl function. It is known
[23], that
\begin{equation}
\label{5.4} m(\lambda)=\frac{1}{\lambda - b_0
-a_0^2\,\widetilde{m}(\lambda)},
\end{equation}
and due to the choice of $a_0,b_0$ (\ref{5.3}) and the definition
of $\widetilde{m}$ we have $$\lambda
-a_0^2\,\widetilde{m}(\lambda) - b_0=\lambda
+\frac{\widetilde{m}(\lambda)}{\alpha A} + \frac{\beta}{\alpha}
+\frac{1}{\alpha \overline{V(i)}}= $$ $$=
\frac{1}{\alpha}\left(f(\lambda) +\frac{1}{
\overline{V(i)}}\right)=\frac{1}{\alpha}\left( \frac{1}{
\overline{V(i)}}-\frac{1}{V(-z(\lambda))}\right).$$ Let $\lambda_k
= \lambda(z_k)=\frac{1}{2}(z_k +z_k^{-1})$, where $z_k$ is taken
from (iv), and so $\lambda_k\in i\mathbb{R}_-$, $\lambda_k\to 0$
as $k\to \infty$. The property (iv) of $V$ yields $\lambda_k - b_0
-a_0^2\,\widetilde{m}(\lambda_k)=0$, that is, $m$ has poles at
$\lambda_k$, and hence $\{\lambda_k\}\in\sigma_d(J)$ and $0$ is a
limit point of $\sigma_d(J)$.

We want to show that there are no other limit points. The argument
rests on the following result of the operator theory.

\begin{lemma}\label{lem5.1}
Let $T$ be a bounded operator in the Hilbert space $\mathcal H$,
such that $T$ and $T^*$ have common cyclic vector $h$. A point
$\lambda\in \mathbb C$ is an isolated point of the spectrum $T$ if
and only if $\lambda$ is a singular point of the Weyl function
$m(\lambda)=(R(\lambda)e_0,e_0)$, $R(\lambda)=(T-\lambda)^{-1}.$
\end{lemma}
{\it Proof.} We have to check that an isolated singularity of the
resolvent $R(\lambda)$ is an unremovable singularity of the Weyl
function(the inverse is obvious). Assume that $m$ is analytic at
$\lambda$. As $$(R(\mu)T^k h,h) = (R(\mu)(T^k - \mu^k)h,h) + \mu^k
m(\mu)=P(\mu) + \mu^k m(\mu),$$ where $P$ is a polynomial, then
the functions $(R(\mu)g_n,h)$ with $g_n=\sum_{k=0}^n c_k T^k h$,
are also analytic at $\lambda$. By the premises the system $\{T^k
h\}_{k\geq 0}$ is complete in $\mathcal{H}$, so that for an
arbitrary $g\in\mathcal{H}$ there is a sequence $g_n\rightarrow g$
as $n\to\infty$. For small enough $\varepsilon>0$ a punctured disk
$\{\mu:\,0<|\mu - \lambda|\leq \varepsilon\}$ lies entirely in the
resolvent set of $T$, i.e., $\|R(\mu)\|\leq C$ for $|\mu -
\lambda|=\varepsilon$. Hence $(R(\mu)g_n,e_0)$ converges to
$(R(\mu)g,e_0)$ uniformly in that disk and so $(R(\mu)g,e_0)$ is
analytic at $\lambda$. Similarly by using the completeness of
$\{(T^*)^k h \}_{k\geq 0}$, we can ascertain that $(R(\mu)g, f)$
is analytic at $\lambda$ for all $g,f \in\mathcal{H}$. It remains
only to note that the weak analyticity implies the strong one
(see, e.g., [24, Chapter V.3]), that is, $R$ is analytic at
$\lambda$, as needed. \hfill $\square$

Assume now that there is a sequence $\{\lambda_k'\}
\in\sigma_d(J)$ and $\lambda'=\lim_{k\to\infty}\lambda_k'\in
[-1,1]\backslash 0$. By lemma \ref{lem5.1} the Weyl function $m$
has poles at $\lambda_k'$ and in view of (\ref{5.3}), (\ref{5.4})
$V(-z(\lambda_k'))=\overline{V(i)}$. But
$-\lim_{k\to\infty}z(\lambda_k') = y'\in\mathbb{T}\setminus \{\pm
i\}$ and since $V$ is analytic at $y'$, then $V\equiv
\overline{V(i)}$, that contradicts  the definition of Pavlov's
function $V$ (\ref{5.1}).

The main problem we have to cope with is to show that $J$ belongs
to $\mathcal{P}(\frac{1}{2} - \varepsilon)$. As a first step in
that direction we find some preliminary rate of stabilization of
the parameters $\{a_n, b_n\}$, which are the three term recurrence
coefficients for orthonormal polynomilas with respect to the
weight $w(x)=A\,\mbox{Im}\,f(x)$ on $[-1,1]$ with $A$ defined in
(\ref{5.21}). As Pavlov's function $V$ is analytic at $1$ and
$V(1)>0$, $V'(1)>0$, then

\begin{equation}\label{5.14}
\mbox{Im}\,V(e^{i\theta}) = V'(1)\sin\theta +
O(\sin^2\theta)\end{equation} for $\theta$ near $0,\,\pi$, and
$$\frac{w(\cos\theta)}{\sin
\theta}\,=\,\frac{\mbox{Im}\,V(e^{i\theta})}{\sin \theta
|V(e^{i\theta})|^2},\quad 0\leq\theta\leq \pi,$$ is $2\pi$
-periodic and infinitely differentiable. Therefore for its Fourier
coefficients
$$\frac{w(\cos\theta)}{\sin \theta} = \sum_{n=-\infty}^\infty\,q_n
e^{i n \theta},\qquad \sum_{n=1}^\infty n^k|q_n|<\infty$$ holds
for any positive integer $k$. According to [25, theorem 1] the
similar inequality is true for
 $\{a_n,\,b_n\}$:
 \begin{equation}\label{5.5}
 \sum_{n=1}^\infty n^k\,\left\{\left|a_n - \frac{1}{2}\right| +
 \left|b_n\right|\right\}\,<\infty,\quad k=1,2,...
 \end{equation}
Relations (\ref{5.5}) make it possible to invoke the scattering
problem technique developed in the previous section. Note that we
do not just pay a tribute to tradition, but apply an adequate
machinery which produces the sharp rate of decay of
$a_n-\frac{1}{2}$ и $b_n$ for weights, related to Pavlov's
function. We do hope that some contemporary methods (such as,
e.g., the Riemann-Hilbert problem) will provide the solution of
the problem we mention for more general classes of weights as
efficient as it is done for the weights analytic in a neighborhood
of $[-1,1]$ \footnote{Our weight $w$ has one singular point at the
origin.} (see [26]).

Let $\{\tilde f_n\}_{n\geq
 0}$ be the Jost solution of the equation $(\tilde J - \lambda)y=0$.
It agrees up to a multiplicative factor which depends on the
spectral parameter, with the Weyl solution for the same equation
and so
 \begin{equation}\label{5.15}\tilde M(z):=-\tilde m(\lambda(z)) =  \frac{\tilde f_1(z)}{\tilde
 f_0(z)}.\end{equation}
By using expression (\ref{4.3}) for the Wronskian $\langle \tilde
f, \overline{\tilde f}\rangle $ we obtain
\begin{equation}\label{5.16} \left|\tilde
f_0\left(e^{i\theta}\right)\right|^2 =
\left|\frac{\sin\theta}{\mbox{Im}\,\tilde
M(e^{i\theta})}\right|.\end{equation}

Let us make sure first that there are no resonance at the edges of
continuous  spectrum of $\tilde J$. Assume on the contrary that
$\tilde f_0(\pm 1)= 0$ at least for one choice of the sign. Since
$A^{-1}\tilde m(\lambda) = - V^{-1}(-z(\lambda)) - \alpha\lambda
-\beta$ is bounded at some neighborhood of $\lambda = \pm 1$, it
follows from (\ref{5.15}) that $\tilde f_1(\pm 1)=0$, and hence by
the recurrence relation (\ref{4.2}) $\tilde f_n(\pm 1)=0$ that
contradicts  the asymptotic behaviour of the Jost solution.

In addition to the properties (i) -- (iv) of Pavlov's function
mentioned above, in [2, formula (3.1)] the bounds for its
derivatives are found
\begin{equation}
\label{5.6} \left|V^{(r)}(e^{i\theta})\right|\leq
K^r\,r!\,r^{\frac{r}{\delta}},\quad \delta = 1 - \gamma,\quad
r=0,1,2,...,\end{equation}
 where $K=K(\gamma)$ is some positive constant. We show that the
 scattering function $\tilde
S = \overline{\tilde f_0}\tilde f_0^{-1}$ obeys the same bounds
with, perhaps, another constant $K$. By [2, lemma 3.П.1 and
remark], when a function $v$ satisfies (\ref{5.6}) and functions
$u$ and $u_1$ are analytic in a neighborhood of the image of $v$
and the closed unit disk $\overline{\mathbb D}$, respectively,
then $u(v)$ and $u_1 v$ satisfy (\ref{5.6}) (with another constant
$K$). That is why the Weyl function
 $$\tilde M(z)= A\left(\frac{1}{V(-z)} +
\frac{\alpha}{2}\left(z+z^{-1}\right) + \beta\right)$$ is subject
to (\ref{5.6}) ($V(e^{i\theta})\neq 0$). Next, since by
(\ref{5.14}) $$\left|\frac{\sin\theta}{\mbox{Im}\tilde
M(e^{i\theta})}\right|\geq C>0$$ for all $\theta$, then by using
(\ref{5.16}) we have
   $$\left|\left[\log\left|\tilde
f_0\left(e^{i\theta}\right)\right|\right]^{(r)}\right|\leq
K_1^r\,r!\,r^{\frac{r}{\delta}},\qquad K_1 =K_1(\gamma)>0,\quad
r\in\mathbb Z_+.$$
 For the scattering function $\tilde S$ the equality
 $$\tilde S(e^{i\theta}) = \exp(-2i\arg\tilde f_0(e^{i\theta})) =
\exp(-2iH(\log|\tilde f_0|)(e^{i\theta})),$$  holds and it remains
to verify that the class (\ref{5.6}) is invariant under the
Hilbert operator $H$.
\begin{lemma}
Let $g$ be a $2\pi$-periodic and infinitely differentiable
function. Then
$$\|\left(H\,g\right)^{(n)}\|_\infty\leq
C\,\left(\|g^{(n)}\|_\infty + \|g^{(n+1)}\|_\infty\right),$$
 where $C$ is a universal constant and the Hilbert operator is
 defined by
\begin{equation}\label{5.7}
(H g)(\theta) = \frac{1}{\pi}\int_0^\pi \frac{g(\theta - t) -
g(\theta + t )}{2\tan\frac{t}{2}}\,d t.
\end{equation}
\end{lemma}
{\it Proof.} In the premises of lemma the integral (\ref{5.7})
converges absolutely and the operator $H$ commutes with the
differentiation $(Hf)^\prime = Hf^\prime$. Next, by the
integration by parts in the partial Fourier sum and taking a limit
we come to the following representation for $f$
$$f(x)=\frac{1}{\pi}\int_{-\pi}^{\pi}f(x) d x +
\frac{1}{\pi}\int_{-\pi}^{\pi} \phi(t)f^\prime(x+t) d t ,$$
 where
$$\phi(t)=\left\{\begin{array}{cc} \frac{t-\pi}{2},\quad
&0<t\leq\pi,\\\frac{t+\pi}{2},\quad &-\pi\leq
t<0.\end{array}\right.$$ Based on this equality we derive $$
\|f\|_\infty\leq C\,\left(\|f\|_2 + \|f^\prime\|_2\right).$$ Since
the Hilbert operator is unitary in $L^2(\mathbb T)$, we have $$\|H
f\|_\infty\leq C\,\left(\|f\|_2 + \|f^\prime\|_2\right)\leq
C\,\left(\|f\|_\infty + \|f^\prime\|_\infty\right),$$
 as was to be proved. \hfill $\square$

Thereby for the scattering function we see that
\begin{equation}\label{5.8}\left|\tilde S^{(r-1)}\left(e^{i\theta}\right)\right|
\leq \,K_2^r\,r!\,r^{\frac{r}{\delta}},\quad
K_2=K_2(\gamma)>0,\quad r=1,2,... .\end{equation}

With lemma \ref{lem4.2} under the belt we can go over to matrix
entries $a_n$, $b_n$. Let us integrate by parts (\ref{4.8}) $r$
times bearing in mind that $\tilde S$ is a periodic function. Then
$$ n^r\,\left(|F(2n+2) - F(2n)| + |F(2n-1) - F(2n+1)|\right)\leq$$
$$\leq
\frac{1}{2\pi\,2^r}\,\int_0^{2\pi}\,\left|\frac{d^r}{d\theta^r}\left[\tilde
S\left(e^{i\theta}\right)\left(e^{2i\theta} -
1\right)\right]\right| d\theta.$$
 The multiplication on a regular function does not alter the structure
of bounds (\ref{5.8}), and so with some positive constant $K_3$
\begin{equation}\label{5.9}|F(2n+2) - F(2n)| + |F(2n-1) - F(2n+1)|\leq C\,
K_3^{r}\,r!\, r^\frac{r}{\delta}\, n^{-r+1}\leq C\, n
\,\inf_{r\geq 3} t^r\, r^{\alpha r}\end{equation}
 with $t=K_3\,n^{-1}\ll 1$ and $\alpha = 1+\frac{1}{\delta}$. We can apply
(\ref{3.20}) to the right hand side of (\ref{5.9}). For large
enough $n$ we have
\begin{equation}\label{5.10}|F(2n+2) - F(2n)| + |F(2n-1) - F(2n+1)|\leq
C \exp \left(-C_1\
n^{\frac{\delta}{1+\delta}}\right),\end{equation} where
 $$C_1 =\frac{1+\delta}{2\delta}\,K_3^{-\frac{\delta}{1+\delta}}
>0.$$
 From (\ref{4.14}) and (\ref{5.10}) it follows that
 $$\hat F^2(2n-2)\leq C\, \exp
\left(-C_1\, n^{\frac{\delta}{1+\delta}}\right)\sum_{k=n}^\infty
\exp \left(-C_1\, k^{\frac{\delta}{1+\delta}}\right)\leq C_2\,
\exp \left(-C_1\, n^{\frac{\delta}{1+\delta}}\right),$$ and by
lemma \ref{lem4.2}
\begin{equation}\label{5.11}
\left|a_n - \frac{1}{2}\right| + |b_n|\leq C_3\, \exp \left(-C_4\,
n^{\frac{\delta}{1+\delta}}\right) = C_3\, \exp \left(-C_4\,
n^{\frac{1-\gamma}{2-\gamma}}\right),\quad C_3, C_4
>0.\end{equation}
 Given $\varepsilon >0$, we choose $\gamma$ so small that
$$\frac{1-\gamma}{2-\gamma}>\frac{1}{2} - \varepsilon,$$
which means exactly $\tilde J\in\mathcal P(1/2-\varepsilon)$.
Theorem \ref{theor2} is now proved for $\nu = 0$.

The general case can be treated  by the same  reasoning by using a
modified Pavlov's function. Indeed, let $-1<\kappa<1$. Consider
the function $$V_\kappa(z)=V\left(\frac{z-\kappa}{1-\kappa
z}\right),\qquad \left(V_0(z) = V(z)\right).$$
 It satisfies (i)--(iv) with the only difference that its singular
 points on the unit circle  $\mathbb T$ are $w_\kappa$ and
 $\overline w_\kappa$, where $w_\kappa =
(i+\kappa)(1+i\kappa)^{-1}$. When $\kappa$ passes through $(-1,1)$
the point $w_\kappa$ describes the upper semicircle. Hence under
an appropriate choice of $\kappa$ we can end up with any point
$\nu,\, -1<\nu<1$ of accumulation for the discrete spectrum of
$J\in\mathcal P(\frac{1}{2} - \varepsilon)$. Now theorem
\ref{theor2} is proved completely.

\setcounter{section}{5} \setcounter{equation}{0}
\setcounter{lemma}{0} \setcounter{theorem}{0}
\section{Appendix.}

We prove here two results concerning perturbation determinants of
complex Jacobi matrices following the line of [6, section 2],
wherein they are established for {\it real symmetric} matrices.

 $1^0.$ {\it Proof of proposition 2.2}. The basic
properties of Schatten--von Neumann operator ideals
$\mathcal{S}_p$ and infinite determinants are presumed to be known
(see, e.g., [5, Chapters III, IV]). The proof is broken up into
steps.

1. Let
$$R_0(z) = (J_0 - \lambda(z))^{-1},\quad
R_0(z)=\{r_{nm}(z)\}_0^\infty,\quad |z|<1, $$ be the resolvent of
the discrete laplacian. Its matrix entries can be found by direct
computation
$$r_{nm}(z)=\frac{z^{|n-m|} - z^{n+m+2}}{z-z^{-1}}
=\,-\sum_{j=0}^{\min(n,m)}\,z^{1+|n-m|+2j},\quad n,m\geq 0,\quad
|z|<1.$$
 It is clear that
\begin{equation}\label{6.1}
|r_{nm}(z)|\leq\min(n,m) +1.
\end{equation}
Take $\zeta\in\mathbb{T}$. The boundary value $r_{nm}(\zeta)$
exists for trivial reasons, however the matrix $\{r_{nm}(\zeta)\}$
is defined only formally (it does not correspond in general to a
bounded operator in $\ell^2$). The situation can be revised by
means of edging $R_0$ with appropriate diagonal matrices. Let
$$D=\mbox{diag} (d_0,d_1,...),\quad d_j\geq 0,\quad
\sum_j\,jd_j<\infty.$$
 Consider an operator function
$$B(z):=D^{1/2}R_0(z)D^{1/2} = \{B_{nm}(z)\}_0^\infty,\quad
|z|<1.$$
 As $D^{1/2}\in\mathcal S_2$ (the Hilbert--Schmidt operator), then
 $B(z)\in\mathcal S_2$ (as a matter of fact,
$B\in\mathcal S_1$ , as a product of two operators from $S_2$, the
fact we will make use of later on). Moreover, by (\ref{6.1})
\begin{equation}\label{6.2}
\|B(z)\|^2_2=\sum_{n,m}\,|B_{nm}(z)|^2\leq\,\sum_{n,m}\,
\left(d_m^{1/2}d_n^{1/2}\,[\min(n,m)+1]\right)^2\leq\,
\left(\sum_j(j+1)d_j\right)^2,\end{equation}
 and hence the boundary values $B(\zeta)=\lim_{z\to\zeta}\,B(z)$
belong to $\mathcal S_2$ and obeys (\ref{6.2}). Finally, $$\|B(z)
- B(\zeta)\|_2^2 = \,\sum_{n,m}\,|B_{nm}(z) -
B_{nm}(\zeta)|^2\,\to 0, \quad z\to\zeta, $$
 by the Lebesgue dominated convergence theorem, and so $B(z)$
 belongs to the disk algebra $\mathcal A$ (as an operator valued
 function in $\mathcal S_2$).

Let now $J$ be the Jacobi matrix (\ref{1.3}), subject to
(\ref{1.4}). Put $d_0 = \max\left(|b_0|,\,\left|c_0 -
\frac{1}{2}\right|\right)$, $ d_n = \max\left(\left|a_{n-1} -
\frac{1}{2}\right|,|b_n|,\left|c_n - \frac{1}{2}\right|\right)$,
 $n\geq 1$ и $D=\mbox{diag}(d_0,d_1,...)$. The perturbation $\Delta J= J-
J_0$ admits factorization $\Delta J = D^{1/2}UD^{1/2}$, where
$$U = \{u_{nm}\},\qquad u_{nm} = \left\{\begin{array}
{cc}\frac{a_{n-1} - \frac{1}{2}}{\sqrt{d_{n-1}d_n}},\quad
&m=n-1,\\ \frac{b_n}{d_n},\quad&m=n,\\\frac{c_{n} -
\frac{1}{2}}{\sqrt{d_{n+1}d_n}},\quad&m=n+1,\end{array}\right.\qquad
u_{nm}=0,\quad |n-m|\geq 2,$$
 where we agree that $\frac{0}{0}=1$.
Obviously, $|u_{nm}|\leq 1$, so that $\|U\|\leq 3$.

2. The function $t(z):=\mbox{tr} (\Delta J R_0(z))$ is well
defined for $|z|<1$ under the condition $\Delta J \in\mathcal
S_1$. Let us show that (\ref{2.1}) yields $t\in \mathcal A$, and
find the upper bound for this function. By the definition of the
trace
 $t=t_1 + t_2 +t_3$, where
$$t_1(z)=\sum_0^\infty\,b_jr_{jj}(z);\
t_2(z)=\sum_0^\infty\,\left(a_j - \frac{1}{2}\right)r_{j,j+1}(z);\
t_3(z)=\sum_0^\infty\,\left(c_j -
\frac{1}{2}\right)r_{j+1,j}(z).$$
 Inequality (\ref{6.1}) along with condition (\ref{2.1})
 enables one to conclude that
$t_l\in\mathcal A,$ $l=1,2,3$ and
\begin{equation}\label{6.3}
|t(z)|\leq\sum_0^\infty\,(j+2)\left\{\left|a_j -
\frac{1}{2}\right| + |b_j| + \left|c_j -
\frac{1}{2}\right|\right\}.\end{equation}

3. Write the perturbation determinant $\Delta$ with factorization
of $\Delta J$ in mind
$$\Delta(z,J) = \det(I + \Delta J R_0)=\det (I +
D^{1/2}UD^{1/2}R_0) = \det (I+UB(z))$$ (we use here the property
$\det(I + AB) = \det (I + BA)$). As $B(z)$ is a nice function
taking the values from $\mathcal S_2$, it looks instructive to go
over to {\it regularized } determinants [5, Chapter IV, \S 2]:
\begin{equation}\label{6.4}
\Delta(z,J) = {\det}_2 (I+UB(z))e^{\mbox{tr}(\Delta J R_0(z))}.
\end{equation}
In view of $\det_2$ being continuous on the operator argument
 [4, theorem IV.2.1] and $t$ being continuous in $
\overline{\mathbb{D}}$, (\ref{6.4}) implies $\Delta\in\mathcal A$.

Apply now (\ref{6.4}) for the associated matrix $J^{(k)}$
$$\Delta(z,J^{(k)})={\det}_2(I+U_kB_k(z))e^{\mbox{tr}(\Delta
J^{(k)}R_0(z))}.$$
 The bordering matrix $D_k$ takes the form $D_k
= \mbox{diag}(d_0^{(k)},d_1^{(k)},...)$, $d_j^{(k)}= d_{j+k}$. By
(\ref{6.2})
$$\|B_k(z)\|^2_2\leq\left(\sum_j\,(j+1)d_j^{(k)}\right)^2$$ and
the Lebesgue theorem provides $\lim_{k\to\infty}\|B_k(z)\|_2 = 0$
uniformly in $\overline{\mathbb{D}}$, whence it follows
$\lim_{k\to\infty}\det_2(I+U_kB_k(z))=1$. By (\ref{6.3})
$$\left|t^{(k)}(z)\right| = \left|\mbox{tr}\left(\Delta
J^{(k)}R_0(z)\right)\right|\leq\,\sum_0^\infty\,(j+2)\left\{\left|a_{j+k}
- \frac{1}{2}\right| + |b_{j+k}| + \left|c_{j+k} -
\frac{1}{2}\right|\right\},$$
 so that $\lim_{k\to\infty}\left|t^{(k)}(z)\right| = 0$ uniformly in
$\overline{\mathbb{D}}$, and the statement is proved.\hfill
$\square$

$2^0.$ {\it Proof of limit relation (\ref{2.2})}. The idea is to
approximate $\Delta J\cdot R_0$ by operators of finite rank and
take into account the continuity of perturbation determinants.

Let $$\hat J_m:=\left(\begin{array}{ccc}J_m &\vdots&
\\\ldots&\ldots&\ldots\\ &\vdots& J_0\end{array}\right),\qquad
\hat J_{0,m}:=\left(\begin{array}{ccc}J_{0,m} &\vdots&
\\\ldots&\ldots&\ldots\\ &\vdots& J_0\end{array}\right)$$
be the block diagonal operators, $$\Delta \hat J_m = \hat J_m -
\hat J_{0,m} = \left(\begin{array}{ccc}J_m - J_{0,m} &\vdots&
\\\ldots&\ldots&\ldots\\ &\vdots& \mathbb{O}\end{array}\right).$$
Since $$\Delta J = J - J_0 = \left(\begin{array}{ccc}J_m
 - J_{0,m} &\vdots&C \\
\cdots&\cdots&\cdots\\A&\vdots&J^{(m)} - J_0\end{array}\right),$$
where the matrices $A$ and $C$ have a single nonzero entry $C_{m1}
= c_m - 1/2$, $A_{1m}=a_m-1/2$, $A_{ji}=C_{ij}=0$, $i\neq m$,
$j\neq 1$, then
$$\|\Delta J - \Delta \hat
J_m\|_1\leq\sum_{j=m}^\infty\,\left\{\left|a_j -
\frac{1}{2}\right| + |b_j| + \left|c_j -
\frac{1}{2}\right|\right\} \to 0, \quad m\to\infty.$$
 Next, it can be easily seen that $J_0 - \hat J_{0,m}$ has only
 two nonzero entries. If $h\in\ell^2$ then
 $$\left(J_0 - \hat J_{0,m}\right) h = (h,e_m)e_{m+1} +
 (h,e_{m+1})e_m \to 0, \quad m\to \infty,$$
or, in other words, $\hat J_{0,m}\rightarrow J_0$ as $ m\to\infty$
in the strong topology (on each vector form $\ell^2$). Let $K$ be
a compact subset of $\mathbb D$. We can match resolvents of
 $J_0$ and $\hat J_{0,m}$ for $z\in K$:
$$\left(\hat J_{0,m} - \lambda\right)^{-1} - \left( J_{0} -
\lambda\right)^{-1}=\left(\hat J_{0,m} -
\lambda\right)^{-1}\left(J_0 - \hat J_{0,m}\right)\left( J_{0} -
\lambda\right)^{-1}.$$
 As $\|\left(\hat J_{0,m} -
\lambda\right)^{-1}\| = (\mbox{dist}(\lambda(z),\sigma(\hat
J_{0,m})))^{-1}$ и $\sigma(\hat J_{0,m}) = [-1,1]$, then
$$\|\left(\hat J_{0,m} - \lambda\right)^{-1}\|\leq C(K), \quad
z\in K.$$
 Hence $\left(\hat J_{0,m} - \lambda\right)^{-1}
\to \left(J_{0} - \lambda\right)^{-1}$ in the strong topology.

According to [5, theorem III.6.3] the convergence of a sequence of
operators $A_m$ to $A$ in $\mathcal S_1$ norm and $X_m$ to $X$ in
the strong topology implies the convergence $A_mX_m$ to $AX$ in
$\mathcal S_1$. Thereby $$\left(\hat J_m - \hat J_{0,m}\right)
\left(\hat J_{0,m} - \lambda\right)^{-1} \to \Delta J R_0,\quad I
+ \left(\hat J_m - \hat J_{0,m}\right) \left(\hat J_{0,m} -
\lambda\right)^{-1} \to I + \Delta J R_0$$ in $\mathcal S_1$ for
each $z$, and hence uniformly on $K$. But the operator $\left(\hat
J_m - \hat J_{0,m}\right)  \left(\hat J_{0,m} -
\lambda\right)^{-1}$ has finite rank and
 $$\det\left(I + \left(\hat J_m - \hat J_{0,m}\right)
 \left(\hat J_{0,m} - \lambda\right)^{-1} \right) =
\det\left(\frac{J_m - \lambda}{J_{0,m} - \lambda}\right),$$ so the
desired result is a consequence of the continuity of perturbation
determinants.

%\newpage

\section{References}

%\bibitem{A} Akhiezer N. Classical moments problem, Fizmatgiz, 1963
%(in Russian).

1. {\it Pavlov B.S. } On nonselfadjoint Schr\"odinger operator I
// In ``Problems of mathematical physics''. LGU, 1966. V. 1.
P.102-132 .

\noindent 2. {\it Pavlov B.S. } On nonselfadjoint Schr\"odinger
operator II
// In ``Problems of mathematical physics''. LGU, 1967. V. 2.
P.135-157 .

\noindent 3. {\it Bairamov E., Cakar O., Krall A.M.}
Non-selfadjoint difference operators and Jacobi matrices with
spectral singularities// Math. Nachr. 2001. V.229. P.5-14.

\noindent 4. {\it Geronimo J.S., Case K.M.} Scattering theory and
polynomials on the real line// Trans. Amer. Math. Soc. 1980.
V.258. P.467-494.

\noindent 5. {\it Gohberg I., Krein M.} Introduction to the theory
of linear nonselfadjoint operators in Hilbert space // М.
``Nauka''. 1965.

\noindent 6. {\it Killip R., Simon B.} Sum rules for Jacobi
matrices and their applications to spectral theory// Ann. Math.
2003. V.158. P.253-321.

\noindent 7. {\it Carleson L.} Sets of uniqueness for functions
analytic in the unit disc// Acta Math. 1952. V.87. P.325-345.

\noindent 8. {\it Aptekarev A.I., Kalyaguine V., Van Assche W.}
Criterion for the resolvent set of nonsymmetric tridiagonal
operators// Proc. Amer. Math. Soc. 1995. V.123. P.2423-2430.

\noindent 9. {\it Beckermann B.} On the convergence of bounded
J--fractions on the resolvent set of the corresponding second
order difference operator // J. Approx. Theory. 1999. V.99.
P.369-408.

\noindent 10. {\it Beckermann B., Kaliaguine V. } The diagonal of
the Pad\'{e} table and the approximation of the Weyl function of
the second order difference operators // Constr. Approx. 1997.
V.13. P.481-510.

\noindent 11. {\it Beckermann B.} Complex Jacobi matrices // J.
Comp. Appl. Math. 2001. V.127. P.17-65.

\noindent 12. {\it Barrios D., Lopez G., Martinez-Finkelstein A.,
Torrano E.} Finite rank approximations of resolvent of infinite
banded matrix and continued fractions // Mat. Sbornik. 1999.
V.190. No. 4. P. 23-42.

\noindent 13. {\it Guseinov G.} Determination of infinite Jacobi
matrix from the scattering data// DAN SSSR. 1976. V.227. No.6.
P.1289 -1292.

\noindent 14. {\it  Teschl G.}   Jacobi Operators and Completely
Integrable Nonlinear Lattices. Math. Surveys and Monographs.
No.72. 1999.

\noindent 15. {\it Fatou P.} S\'{e}ries trigonom\'{e}triques et
s\'{e}ries de Taylor // Acta Math. 1906. V.30. P.335-400.

\noindent 16. {\it Beurling A.} Ensembles exceptionneles// Acta
Math. 1939. V.72. P.1-13.

\noindent 17. {\it Khrushchev S.} The problem of simultaneous
approximation and erasure of singularities for Cauchy type
integrals// Trudy MIAN SSSR. 1978. V.130. P.124-195.

\noindent 18. {\it Besicovitch A.S., Taylor S.J.} On the
complementary intervals of a linear closed sets of zero Lebesque
measure // J. Lond. Math. Soc. 1954. V.29. P.449-459.

\noindent 19. {\it Taylor B.A., Williams D.L.} Boundary zero sets
of $A^\infty$ functions satisfying growth conditions // Proc.
Amer. Math. Soc. 1972. V.35. P.155-160.

\noindent 20. {\it Korenblum B.} Quasi-analytic classes of
functions in a disk // DAN SSSR 1965. V.164. P.36-39.

\noindent 21. {\it Nikishin Е. } The discrete Sturm--Liouville
operatr and some probof the function theory // Trudy Petrovsky
seminar. 1984. V. 10. P.3-77.

\noindent 22. {\it  Toda М. } Theory of nonlinear lattices.
 Moscow: "Mir". 1984.

\noindent 23. {\it Kaliagin V.} On rational approximation of the
resolvent function of difference second order operator // UMN.
1994. V.49. No. 3. P.181-182.

\noindent 24. {\it Yosida K. } Functional analysis//
Springer-Verlag. 1965.

\noindent 25. {\it Geronimo J.S. } A relation between the
coefficients in the recurrence formula and the spectral function
for orthogonal polynomials // Trans. Amer. Math. Soc. 1980. V.260.
No.1. P.65-82.

\noindent 26. {\it Kuijlaars A., McLaughlin  K.T.-R., Van Assche
W., Vanlessen  M.} The Riemann-Hilbert approach to strong
asymptotics for orthogonal polynomials on $[-1, 1]$
// Adv. in Math. 2004. V. 188. P.337-398.

\end{document}